\documentclass[a4paper,12pt]{amsart}
\usepackage{amssymb}
\usepackage{ifthen}
 \usepackage[dvips]{graphicx}
\nonstopmode \numberwithin{equation}{section}
\usepackage{amssymb}
\usepackage{ifthen}
\usepackage{graphicx}
\usepackage{amsmath}
\usepackage[T1]{fontenc} %skandit
\usepackage[utf8]{inputenc}
\usepackage[usenames,dvipsnames]{color}
\usepackage{color}
\usepackage[english]{babel}
\usepackage{fancyhdr}
\usepackage{fancybox}
\usepackage{tikz}
\nonstopmode \numberwithin{equation}{section}
\theoremstyle{plain}
\newtheorem{thm}[equation]{Theorem}
\newtheorem{cor}[equation]{Corollary}
\newtheorem{lem}[equation]{Lemma}
\newtheorem{prop}{Proposition}

\newtheorem{conj}{Conjecture}

\theoremstyle{definition}
\newtheorem{defn}{Definition}[section]

\newtheorem{prob}{Problem}
\newtheorem{rem}{Remark}[section]

\newcounter{minutes}
\divide\time by 60
\newcounter{hours}
\multiply\time by 60
\addtocounter{minutes}{-\time}
\newcounter {own}
\def\theown {\thesection       .\arabic{own}}

\newenvironment{pf}[1][]{%
 \vskip 3mm
 \noindent
 \ifthenelse{\equal{#1}{}}%
  {{\slshape Proof. }}%
  {{\slshape #1.} }%
 }%
{\qed\bigskip}

\newcounter{alphabet}

\def\be{\begin{equation}}
\def\ee{\end{equation}}

\newcommand{\bee}{\begin{enumerate}}
\newcommand{\eee}{\end{enumerate}}

\newcommand{\blem}{\begin{lem}}
\newcommand{\elem}{\end{lem}}
\newcommand{\bthm}{\begin{thm}}
\newcommand{\ethm}{\end{thm}}
\newcommand{\bcor}{\begin{cor}}
\newcommand{\ecor}{\end{cor}}
\newcommand{\beg}{\begin{examp}}
\newcommand{\eeg}{\end{examp}}
\newcommand{\begs}{\begin{examples}}
\newcommand{\eegs}{\end{examples}}

\newcommand{\bdefn}{\begin{defn}}
\newcommand{\edefn}{\end{defn}}

\newcommand{\bprob}{\begin{prob}}
\newcommand{\eprob}{\end{prob}}
\newcommand{\bei}{\begin{itemize}}
\newcommand{\eei}{\end{itemize}}

\newcommand{\bcon}{\begin{conj}}
\newcommand{\econ}{\end{conj}}
\newcommand{\bcons}{\begin{conjs}}
\newcommand{\econs}{\end{conjs}}
\newcommand{\bprop}{\begin{prop}}
\newcommand{\eprop}{\end{prop}}
\newcommand{\br}{\begin{rem}}
\newcommand{\er}{\end{rem}}
\newcommand{\brs}{\begin{rems}}
\newcommand{\ers}{\end{rems}}
\newcommand{\bo}{\begin{obser}}
\newcommand{\eo}{\end{obser}}
\newcommand{\bos}{\begin{obsers}}
\newcommand{\eos}{\end{obsers}}
\newcommand{\bpf}{\begin{pf}}
\newcommand{\epf}{\end{pf}}
\newcommand{\ba}{\begin{array}}
\newcommand{\ea}{\end{array}}
\newcommand{\beq}{\begin{eqnarray}}
\newcommand{\beqq}{\begin{eqnarray*}}
\newcommand{\eeq}{\end{eqnarray}}
\newcommand{\eeqq}{\end{eqnarray*}}

\begin{document}

\title{Bohr inequalities for unimodular bounded functions on simply connected domains}

\author{Molla Basir Ahamed}
\address{Molla Basir Ahamed,
	School of Basic Science,
	Indian Institute of Technology Bhubaneswar,
	Bhubaneswar-752050, Odisha, India.}
\email{mba15@iitbbs.ac.in}

\author{Vasudevarao Allu}
\address{Vasudevarao Allu,
School of Basic Science,
Indian Institute of Technology Bhubaneswar,
Bhubaneswar-752050, Odisha, India.}
\email{avrao@iitbbs.ac.in}

\author{Himadri Halder}
\address{Himadri Halder,
School of Basic Science,
Indian Institute of Technology Bhubaneswar,
Bhubaneswar-752050, Odisha, India.}
\email{hh11@iitbbs.ac.in}

\subjclass[{AMS} Subject Classification:]{Primary 30C45, 30C50, 30C80}
\keywords{Unimodular bounded functions, analytic functions, harmonic mappings; Simply connected domains, improved Bohr inequality, Bohr radius.}

\def\thefootnote{}
\footnotetext{ {\tiny File:~\jobname.tex,
printed: \number\year-\number\month-\number\day,
          \thehours.\ifnum\theminutes<10{0}\fi\theminutes }
} \makeatletter\def\thefootnote{\@arabic\c@footnote}\makeatother

\begin{abstract}
	Let $ \mathcal{H}(\mathbb{D}) $ be the class of analytic functions in the unit disk $ \mathbb{D} : =\{z\in\mathbb{C} : |z|<1\} $. The classical Bohr's inequality states that if a power series $ f(z)=\sum_{n=0}^{\infty}a_nz^n $ converges in $ \mathbb{D} $ and $ |f(z)|<1 $ for $ z\in\mathbb{D} $, then 
	\begin{equation*}
		\sum_{n=0}^{\infty}|a_n|r^n\leq 1\;\;\mbox{for}\;\; r\leq \frac{1}{3}
	\end{equation*}
	and the constant $ 1/3 $ cannot be improved.  The constant $ 1/3 $ is known as Bohr radius. In this paper, we study Bohr phenomenon for analytic as well as harmonic mappings on simply connected domains. We prove several sharp results on improved Bohr radius for analytic functions as well as for harmonic mappings on simply connected domains. 
\end{abstract}

\maketitle
\pagestyle{myheadings}
\markboth{Molla Basir Ahamed, Vasudevarao Allu and  Himadri Halder}{Bohr inequalities for unimodular bounded functions on simply connected domains}

\section{Introduction}
\vspace{2mm}
Let $ \Omega\subseteq\mathbb{C} $ be a simply connected domain containing the unit disk $\mathbb{D}:=\{z\in \mathbb{C}: |z|<1\}$ and $ \mathcal{B}(\Omega)$ be the class of analytic functions $f$ in $\Omega$ such that $|f(z)|<1$ for all $z \in \Omega$. The Bohr radius $ B_{\Omega} $ for the class $ \mathcal{B}(\Omega) $ is defined by
\begin{equation*}
B_{\Omega}:=\sup\bigg\{r\in (0,1) : \sum_{n=0}^{\infty}|a_n|r^n\leq 1\; \text{for all}\; f(z)=\sum_{n=0}^{\infty}a_nz^n\in\mathcal{B}(\Omega),\; z\in\mathbb{D}\bigg\}.
\end{equation*}
For $f \in  B_{\Omega}$ given by $f(z)=\sum_{n=0}^{\infty}a_nz^n$ in $\mathbb{D}$, its associated majorant series defined by $M_{f}(r):=\sum_{n=0}^{\infty}|a_n|r^n$.
For $\Omega = \mathbb{D}$, it is proved that $ B_{\mathbb{D}}=1/3 $, which is the classical Bohr radius. The classical result of Bohr \cite{Bohr-1914}, which inspired a lot in the recent years is in the following form after subsequent improvements due to M. Riesz, I. Schur, and N. Wiener.
 \begin{thm}\cite{Bohr-1914} \label{thm-1.1}
 	If $ f\in\mathcal{B}(\mathbb{D}) $ and $ f(z)=\sum_{n=0}^{\infty}a_nz^n $ for $z \in \mathbb{D}$, then 
 	 \begin{equation}\label{e-1.2a}
 	 	\sum_{n=0}^{\infty}|a_n|r^n\leq 1\;\;\mbox{for}\;\; |z|=r\leq \frac{1}{3}.
 	 \end{equation}
  The radius $ 1/3 $ cannot be improved.
 \end{thm}
It is important to note that if $|f(z)|\leq 1$ in $\mathbb{D}$ and $|f(z^*)|=1$ for some $z^*$ in $\mathbb{D}$, then $f(z)$ reduces to a unimodular constant. Thus, we restrict our attention to $|f(z)|<1$. It is worth to point out that if $|a_0|$ in the equality \eqref{e-1.2a} is replaced by $|a_0|^2$, then the radius $1/3$ could be replaced by $1/2$. Furthermore, if $a_0=0$ in Theorem \ref{thm-1.1}, then the sharp radius can be improved to $1/\sqrt{2}$. In a more convenient way, we can demonstrate these facts that a minor change in the coefficients gives different sharp radius, and thus, sharp coefficient estimates which play a vital role to obtain the sharp radius. Therefore, proofs of these results and the inequality \eqref{e-1.2a} relied on the sharp coefficient inequalities which may be obtained as an application of Pick's invariant form of Schwarz's lemma for $ f\in\mathcal{H}(\mathbb{D}) $ 
 \begin{equation*}
 	|f^{\prime}(z)|\leq\frac{1-|f(z)|^2}{1-|z|^2}\;\;\mbox{for}\;\; z\in\mathbb{D}.
 \end{equation*}  
 In particular, $ |f^{\prime}(0)|=|a_1|\leq 1-|f(0)|^2=1-|a_0|^2 $ and hence from this, the sharp inequality $ |a_n|\leq 1-|a_0|^2 $ follows for $ n\geq 1 $.
\\

 Studying the Bohr inequalities in recent years become an interesting topic of research in functions of one as well as several complex variables. The notion of Bohr inequality has been generalized to several complex variables, to planar harmonic mappings, to polynomials, to solutions of elliptic partial differential equations, and more abstract settings.  In $1997$, Boas and Khavinson \cite{boas-1997} extended the Bohr inequality \eqref{e-1.2a} to several complex variables by finding multidimensional Bohr radius. For more interesting aspects of Bohr inequalities in this particular direction, we refer the reader to \cite{aizenberg-2001,aizn-2007,bene-2004,boas-1997}. For various forms of Bohr inequalities we refer to (see \cite{Huang-Liu-Ponnu, Ismagilov-2020}). Many interesting multidimensional analogue of improved Bohr inequality have been obtained by Liu {\it et al.} (see \cite{Liu-Pon-PAMS-2020}). Recently, Bohr inequalities for functions defined in a simply connected domain have been developed by Evdoridis {\it et al.} \cite{Evd-Ponn-Rasi-2020} (also see \cite{Ahamed-Allu-Halder-P3-2020}). For various research works on Bohr inequality as stated above, we suggest the reader glance through the articles \cite{aizn-2000, alkhaleefah-2019, Himadri-Vasu-P2} and the references therein.
 \vspace{3mm} 
 
 Let $\mathbb{D}(a,r):=\{z \in \mathbb{C}: |z-a|<1\}$. For $0\leq \gamma <1$, we consider the disk $\Omega _{\gamma}$ defined by 
$$
\Omega_{\gamma}:=\bigg\{z\in\mathbb{C} : \bigg|z+\frac{\gamma}{1-\gamma}\bigg|<\frac{1}{1-\gamma}\bigg\}.
$$
It is easy to see that $\Omega _{\gamma}$ always contains the unit disk $\mathbb{D}$. The following coefficient estimates for class $\mathcal{B}(\Omega_{\gamma})$ are required to prove some of our results.
\begin{lem}\cite{Evd-Ponn-Rasi-2020}\label{lem-2.1}
    For $ \gamma\in[0,1) $, let 
	\begin{equation*}
	\Omega_{\gamma}:=\bigg\{z\in\mathbb{C} : \bigg|z+\frac{\gamma}{1-\gamma}\bigg|<\frac{1}{1-\gamma}\bigg\},
	\end{equation*}
	and let $ f $ be an analytic function in $ \Omega_{\gamma} $, bounded by $ 1 $, with the series representation $ f(z)=\sum_{n=0}^{\infty}a_nz^n $ in the unit disk $ \mathbb{D}. $ Then 
	\begin{equation*}
	|a_n|\leq\frac{1-|a_0|^2}{1+\gamma}\quad\mbox{for}\quad n\geq 1.
	\end{equation*}
\end{lem}
In $ 2010 $, the notion of classical Bohr inequality $ \sum_{n=0}^{\infty}|a_n|r^n\leq 1 $ has been generalized by Fournier and Ruscheweyh \cite{Four-Rusc-2010} to the class $\mathcal{B}(\Omega_{\gamma})$. More precisely,
\begin{thm}\cite{Four-Rusc-2010} \label{thm-1.2}
	For $ 0\leq \gamma<1 $, let $ f\in\mathcal{B}(\Omega_{\gamma}) $, with $ f(z)=\sum_{n=0}^{\infty}a_nz^n $ in $ \mathbb{D} $. Then,
	\begin{equation*}
		\sum_{n=0}^{\infty}|a_n|r^n\leq 1\;\; \text{for}\;\; r\leq\rho:=\frac{1+\gamma}{3+\gamma}.
	\end{equation*}
	Moreover, $ \sum_{n=0}^{\infty}|a_n|\rho^n=1 $ holds for a function $ f(z)=\sum_{n=0}^{\infty}a_nz^n $ in $ \mathcal{B}(\Omega_{\gamma}) $ if, and only if, $ f(z)=c $ with $ |c|=1 $.
\end{thm}
In $ 2010 $, Fournier and Ruscheweyh \cite{Four-Rusc-2010} proved the following result.
\begin{thm}\cite{Four-Rusc-2010}
	Let $ \Omega $ be a simply connected domain which contains the unit disk $ \mathbb{D} $ and let 
	\begin{equation}\label{e-2.3}
		\lambda:=\lambda(\Omega)=\sup_{f\in\mathcal{B}(\Omega),\;n\geq 1}\bigg\{\frac{|a_n|}{1-|a_0|^2} :  a_0\neq f(z)=\sum_{n=0}^{\infty}a_nz^n\;\; \mbox{for}\;\; z\in\mathbb{D}\bigg\}.
	\end{equation}
	Then $ 1/(1+2\lambda)\leq\mathcal{B}_{\Omega} $ and the equality $ \sum_{n=0}^{\infty}|a_n|(1/(1+2\lambda))^n=1 $ holds for a function $ f(z)=\sum_{n=0}^{\infty}a_nz^n $ in $ \mathcal{B}(\Omega) $ if, and only if, $ f\equiv c $ with $ |c|=1 $.
\end{thm}
A complex-valued function $ f $ in a simply connected domain $ \Omega\subseteq\mathbb{C} $ is called a harmonic in $ \Omega $ if it satisfies the Laplace equation $ \Delta f=4f_{z\overline{z}}=0 $. It is known that each harmonic mapping $ f $ has a canonical representation of the form $ f=h+\overline{g} $ where $ h $ and $ g $ are analytic in a simply connected domain $ \Omega\subseteq\mathbb{C} $ and this representation is unique up to an additive constant. A locally univalent function $ f $ is sens-preserving if its Jacobian $ J_f(z)>0 $ in $ \Omega $, where $ J_f $ is defined by $ J_f(z)=|h^{\prime}(z)|^2-|g^{\prime}(z)|^2 $.  Every harmonic function $ f $ in $ \mathbb{D} $ has the following representation
\begin{equation}\label{e-2.4a}
f(z)=h(z)+\overline{g(z)}=\sum_{n=0}^{\infty}a_nz^n+\overline{\sum_{n=0}^{\infty}b_nz^n}\;\;\mbox{for}\;\; z\in\mathbb{D}.
\end{equation} 
In the recent years, harmonic extensions of the classical Bohr radius have been extensively studied by several authors and we refer the reader to glance through the articles \cite{Evdoridis-Ponnusamy-2018,Kayumov-Ponnusamy-2018-b,Kayumov-Ponnuswamy-MN-2018,Liu-Ponnusamy-BMMS-2019,Liu-Ponnusamy-Wang-2020}. 
%In $ 2020 $, Evdoridis \textit{et al.} \cite{Evd-Ponn-Rasi-2020} proved the following interesting lemma for simply connected domain $ \Omega_{\gamma} $.

\section{Main Results}
\noindent Before we state an improved version of Theorem \ref{thm-1.2} for class $\mathcal{B}(\Omega_{\gamma})$, we prove the following lemma. Let $ Q(w) $ be a polynomial defined by
\begin{equation}\label{e-1.8}
Q(w)=c_1w+c_2w^2+\ldots+c_mw^m\;\; \mbox{for}\;\; c_j\in\mathbb{R}^{+}\;\; j=1, 2, \cdots, m.
\end{equation}

\begin{lem}\label{lem-3.6}
	Let $ g : \mathbb{D}\rightarrow\overline{\mathbb{D}} $ be an analytic function and $ \gamma\in\mathbb{D} $ be such that $ g(z)=\sum_{n=0}^{\infty}\alpha_n(z-\gamma)^n $ for $ |z-\gamma|<1-|\gamma| $. Then
	\begin{align*}
	\sum_{n=0}^{\infty}|\alpha_n|\rho^n+Q\left(\frac{S_{\rho}}{\pi}\right)\leq 1\;\; \text{for}\;\; \rho\leq \rho_0=\frac{1-|\gamma|^2}{3+|\gamma|},
	\end{align*}
	where $ S_{\rho} $ denotes the area of the image of the disk $ \mathbb{D}(\gamma;\rho(1-|\gamma|)) $ under the mapping  $ g $ and the non-negative real coefficients $ c_1, c_2, \cdots, c_m $ of the polynomial $ Q(w) $ given by \eqref{e-1.8} satisfy
	\begin{equation*}
	8c_1\left(\frac{3}{8}\right)^2+24c_2\left(\frac{3}{8}\right)^4+\cdots+8(2m-1)c_m\left(\frac{3}{8}\right)^{2m}=1.
	\end{equation*}
\end{lem}

\noindent In view of Lemma \ref{lem-3.6}, we obtain the following  improved version of Theorem \ref{thm-1.2}. In deed, with an additional polynomial $ Q\left(S_{r(\gamma-1)}/\pi\right) $ with the majorant series of $ f $, we prove that the sum is still less than or equals to $1$ for the radius $r_{0}=(1+\gamma) /(3+\gamma) $. More precisely,  
\begin{thm}\label{th-3.6}
	For $ 0\leq\gamma<1, $ let $ f\in\mathcal{B}(\Omega_{\gamma}) $ with $ f(z)=\sum_{n=0}^{\infty}a_nz^n $ in $ \mathbb{D} $, then we have 
	\begin{equation*}
	\sum_{n=0}^{\infty}|a_n|r^n+Q\left(\frac{S_{r(\gamma-1)}}{\pi}\right)\leq 1\;\;\mbox{for}\;\; r\leq r_0=\frac{1+\gamma}{3+\gamma},
	\end{equation*} where the coefficients of $ Q(w) $ connected by the relations
	\begin{equation*}
	8c_1\left(\frac{3}{8}\right)^2+24c_2\left(\frac{3}{8}\right)^4+\cdots+8(2m-1)c_m\left(\frac{3}{8}\right)^{2m}=1.
	\end{equation*} Furthermore, the quantities $ c_1, c_2, \cdots,c_m $ and $ (1+\gamma)/(3+\gamma) $ can not be improved.
\end{thm}

Let $ P(w) $ be a polynomial of degree $ m $ defined by
\begin{equation}\label{e-2.3a}
	P(w)=k_1w+k_2w^2+\cdots+k_mw^m,\;\;\mbox{where}\;\; k_j=\left(\frac{1+\lambda}{1+2\lambda}\right)^{2j}\;
\end{equation}  
for $ j=1, 2, \cdots, m. $ Our next aim is to prove a more general version of Theorem \ref{th-3.6}, where an improved Bohr inequality is obtained for the analytic functions defined on a simply connected domain $\Omega\subseteq\mathbb{C}$.
\begin{thm}\label{th-3.3}	
Let $ \Omega $ be a simply connected domain containing the unit disk $ \mathbb{D} $ and $ f\in\mathcal{B}(\Omega) $ with $ f(z)=\displaystyle\sum_{n=0}^{\infty}a_nz^n $ for $ z\in\mathbb{D} $. Then for the polynomial $ P(w) $ 
	\begin{equation*}
		\sum_{n=0}^{\infty}|a_n|r^n+P\left(\frac{S_r}{\pi}\right)\leq 1\;\; \mbox{for}\;\; r\leq\frac{1}{1+2\lambda},
	\end{equation*}
where $ S_r $ denotes the area of the image of the disk $ \mathbb{D}(0,r) $ under the mapping $ f $ and $ P(w) $ is a polynomial given by \eqref{e-2.3a}. The equality $ \sum_{n=0}^{\infty}|a_n|r^n+P\left({S_r}/{\pi}\right)=1 $ holds for the function $ f\in \mathcal{B}(\Omega) $ if, and only if, $f$ is an unimodular function {\it i.e.}$, f\equiv c $ with $ |c|=1. $
\end{thm}

Next, we show that if we replace $|a_n|$ by $|a_n|+\beta |a_n|^2$ in Theorem \ref{th-3.3}, then we obtain the same radius $1/(1+2\lambda)$ as we obtained in Theorem \ref{th-3.3}.  	
\begin{thm}\label{th-3.4}
	Let $ \Omega $ be a simply connected domain containing the unit disk $ \mathbb{D} $ and $ f\in\mathcal{B}(\Omega) $ with $ f(z)=\sum_{n=0}^{\infty}a_nz^n $ for $ z\in\mathbb{D} $. Then for $ 0\leq\beta\leq 1/\lambda $,
	\begin{equation*}
		|a_0|+\sum_{n=1}^{\infty}\bigg(|a_n|+\beta |a_n|^2\bigg)r^n\leq 1\;\; \mbox{for}\;\; r\leq\frac{1}{1+2\lambda}.
	\end{equation*}
The equality $ |a_0|+\sum_{n=1}^{\infty}\bigg(|a_n|+\beta|a_n|^2\bigg)(1/(1+2\lambda))^n=1 $ holds for a function $ f\in \mathcal{B}(\Omega) $ if, and only if, $ f\equiv c $ with $ |c|=1 $.
\end{thm}
By refining the coefficients of functions in the class $ \mathcal{B}(\Omega) $, we obtain the following interesting result.
\begin{thm}\label{th-3.5}
	Let $ \Omega $ be a simply connected domain containing the unit disk $ \mathbb{D} $ and $ f\in\mathcal{B}(\Omega) $ with $ f(z)=\sum_{n=0}^{\infty}a_nz^n $ for $ z\in\mathbb{D} $. Then 
	\begin{equation*}
		\sum_{n=0}^{\infty}|a_n|r^n+\left(\frac{1+\lambda}{2\lambda(1+|a_0|)}+\frac{2}{3}\frac{(1+\lambda)r}{1-r}\right)\sum_{n=1}^{\infty}|a_n|^2r^{2n}\leq 1\;\; \mbox{for}\;\; r\leq\frac{1}{1+2\lambda}.
	\end{equation*}
	The equality 
	\begin{equation*}
		\sum_{n=0}^{\infty}|a_n|\left(\frac{1}{1+2\lambda}\right)^n+\left(\frac{1+\lambda}{2\lambda(1+|a_0|)}+\frac{2}{3}\frac{(1+\lambda)}{2\lambda}\right)\sum_{n=1}^{\infty}|a_n|^2\left(\frac{1}{1+2\lambda}\right)^{2n}=1
	\end{equation*}
	 holds for a function $ f\in\mathcal{B}(\Omega) $ if, and only if, $ f$ is an unimodular function.
\end{thm}

\begin{thm}\label{th-3.1}
	Let $ f=h+\overline{g} $ given by \eqref{e-2.4a} be harmonic in $ \Omega_{\gamma} $, with $ |h(z)|+|g(z)|\leq 1$ in $ \Omega_{\gamma} $. Then 
	\begin{equation*}
		|a_0|+\sum_{n=1}^{\infty}\bigg(|a_n|+|b_n|\bigg)r^n\leq 1 \quad\mbox{for}\quad |z|=r\leq r_0 : = \frac{1+\gamma}{3+\gamma}.
	\end{equation*}
\end{thm}
\begin{thm}\label{th-3.2}
	Let $ f=h+\overline{g} $ given by \eqref{e-2.4a} be harmonic in $ \Omega_{\gamma} $, with $ |h(z)|+|g(z)|\leq 1$ in $ \Omega_{\gamma} $. Then 
	\begin{equation*}
		\sqrt{|a_0|^2+|b_0|^2}+\sum_{n=1}^{\infty}\sqrt{|a_n|^2+|b_n|^2}r^n\leq 1 \quad\mbox{for}\quad |z|=r\leq r_0 := \frac{1+\gamma}{3+\gamma}.
	\end{equation*}
\end{thm}

\section{Proof of the main results}
\begin{pf}[\bf Proof of Lemma \ref{lem-3.6}]
	Without loss of generality, we assume that $ \gamma\in [0,1) $. Then it is easy to see that $ z\in\mathbb{D}_{\gamma}:=\mathbb{D}(\gamma;1-\gamma) $ if, and only if, $ w=(z-\gamma)/(1-\gamma)\in\mathbb{D}. $
	Then we have 
	\begin{align*}
		g(z)=\sum_{n=0}^{\infty}\alpha_n(1-\gamma)^n\phi^n(z)=\sum_{n=0}^{\infty}b_n\phi^n(z):=G(\phi(z))
	\end{align*}
	for $ z\in \mathbb{D}_{\gamma}, $ where $ b_n= \alpha_{n} (1-\gamma)^n$. A simple computation shows that
	\begin{align}\label{e-2.8}
		\frac{S_{\rho}}{\pi}=\frac{1}{\pi}\text{Area}\bigg(G(\mathbb{D}(0,\rho))\bigg)\leq (1-|b_0|^2)^2\frac{\rho^2}{(1-\rho^2)^2}=(1-|\alpha_0|^2)^2\frac{\rho^2}{(1-\rho^2)^2}.
	\end{align}
	Therefore,
	\begin{align}\label{e-2.9}
		\sum_{n=1}^{\infty}|\alpha_n|\rho^n\leq\frac{1-|\alpha_0|^2}{1+\gamma}\sum_{n=1}^{\infty}\left(\frac{\rho}{1-\gamma}\right)^n=\frac{1-|\alpha_0|^2}{1+\gamma}\left(\frac{\rho}{1-\gamma-\rho}\right).
	\end{align}
	In view of \eqref{e-2.8} and \eqref{e-2.9}, we obtain 
	\begin{align*} &
		|\alpha_0|+	\sum_{n=1}^{\infty}|\alpha_n|\rho^n+Q\left(\frac{S_{\rho}}{\pi}\right)\\&= |\alpha_0|+\frac{(1-|\alpha_0|^2)\rho}{(1+\gamma)(1-\gamma-\rho)}+\sum_{j=0}^{m}c_j\bigg(\frac{(1-|\alpha_0|^2)\rho}{(1-\rho^2)}\bigg)^{2j}\\&= 1+\Psi_1^{\gamma}(\rho),
	\end{align*}
	where 
	\begin{align*}
		\Psi_1^{\gamma}(\rho)=\frac{(1-|\alpha_0|^2)\rho}{(1+\gamma)(1-\gamma-\rho)}+\sum_{j=0}^{m}c_j\bigg(\frac{(1-|\alpha_0|^2)\rho}{(1-\rho^2)}\bigg)^{2j}-(1-|\alpha_0|).
	\end{align*}
	Let $ G_m(\rho) $ be defined by
	\begin{align*}
		G_m(\rho)&=\left(\frac{c_{m-1}}{c_m}\right)\frac{\rho^{2m-2}}{(1-\rho^2)^{2m-2}(1-|\alpha_0|^2)^{2m-2}}+\cdots+\left(\frac{c_{m-j}}{c_m}\right)\frac{\rho^{2m-2j}}{(1-\rho^2)^{2m-2j}(1-|\alpha_0|^2)^{2j}} \\[2mm]
		&\quad\quad+\cdots+\left(\frac{c_1}{c_m}\right)\frac{\rho^2}{(1-\rho^2)^2(1-|\alpha_0|^2)^{2m-2}}+\frac{\rho^{2m}}{(1-\rho^2)^{2m}}.
	\end{align*}
	Then we can write $ \Psi_1^{\gamma}(\rho) $ as
	\begin{align*} 
		\Psi_1^{\gamma}(\rho)&=\frac{1-|\alpha_0|^2}{2}\bigg(1+2c_m(1-|\alpha_0|^2)^{2m-2}G_m(\rho)-\frac{2}{1+|\alpha_0|} \\&\quad\quad+\frac{1}{2\lambda(1-|\alpha_0|^2)^3}\left(\frac{2\rho}{(1+\gamma)(1-\gamma-\rho)}-1\right)\bigg).
	\end{align*}
	We suppose that $ \rho\leq \rho_0 $. Then it is easy to see that $\Psi_1^{\gamma}(\rho)$ is an increasing function and hence $ 	\Psi_1^{\gamma}(\rho)\leq 	\Psi_1^{\gamma}(\rho_0) $, where 
	\begin{align*}
		\frac{2\rho_0}{(1+\gamma)(1-\gamma-\rho_0)}=1
	\end{align*}
	which is equivalent to 
	\begin{equation*}
	\rho_0=\frac{1-\gamma^2}{3+\gamma}.
	\end{equation*}
	A simple computation shows that 
	\begin{align*}
		\Psi_1^{\gamma}(\rho_0)=\frac{1-|\alpha_0|^2}{2}\left(1+2F_m(|\alpha_0|)-\frac{2}{1+|\alpha_0|}\right):=\frac{1-|\alpha_0|^2}{2}J(|\alpha_0|),
	\end{align*}
	where 
	\begin{align*}
		F_m(|\alpha_0|)&=c_m(1-|\alpha_0|^2)^{2m-1}A^{2m}(\gamma)+c_{m-1}(1-|\alpha_0|^2)^{2m-3}A^{2m-1}(\gamma)+\cdots\\&\quad\quad+c_{m-j}(1-|\alpha_0|^2)^{2m-2j-1}A^{2j-1}(\gamma)+\cdots+c_1(1-|\alpha_0|^2)A^{2}(\gamma),\\
		J(x)&=1+2F_m(x)-\frac{2}{1+x}\;\; \mbox{for}\;\ x\in [0,1]\;\;\text{and}\;\;\\ A(\gamma)&=\frac{(3+\gamma)(1-\gamma^2)}{(3+\gamma)^2-(1-\gamma^2)^2}.
	\end{align*}
	 We note that $ A(\gamma)>0 $ for $ \gamma\in [0,1) $. In order to show that $ \Psi_1^{\gamma}(\rho_0)\leq 0 $, it is enough to show that $ J(x)\leq 0 $ for $ x\in [0,1] $ and $ \gamma\in [0,1) $ so that $ 	\Psi_1^{\gamma}(\rho_0)\leq 0 $. Further, a simple computation shows that
	\begin{align*}
		J(0)&=2c_mA^{2m}(\gamma)+2c_{m-1}A^{2m-1}(\gamma)+\cdots+c_{m-j}A^{2j-1}(\gamma)+\cdots+c_1A^2(\gamma)-1\\ \;\; \text{and}\;\; &\lim_{x\rightarrow 1^{-}}J(x)=0.
	\end{align*}
	It is easy to see that $ A(\gamma)=(f_1\circ f_2)(\gamma) $, where $ f_1(\rho)=\rho/(1-\rho^2) $ and $ f_2(\gamma)=(1-\gamma^2)/(3+\gamma). $
	Since $ A^{\prime}(\gamma)=f^{\prime}_1(f_2(\gamma))f^{\prime}_2(\gamma) $ and 
	\begin{equation}
		f^{\prime}_2(\gamma)=-\left(\frac{\gamma^2+6\gamma+1}{(3+\gamma)^2}\right)<0,
	\end{equation}
	we show that $ f_1(\rho) $ is an increasing function of $ \rho $ in $ (0,1) $, and $ f_2 $ is a decreasing function of $ \gamma $ in $ [0,1) $ so that $ A(\gamma) $ is a decreasing function of $ \gamma $  in $ [0,1) $, with $ A(0)=3/8 $ and $ A(1)=0 $. It can be seen that each $ A^{2j}(\gamma) $ is a decreasing function on $ [0,1) $, for $ j=1, 2, \cdots, m $. Therefore, we have 
	\begin{align*}
		A^{2j}(\gamma)\leq A^{2j}(0)=\left(\frac{3}{8}\right)^{2j}\;\; \mbox{for}\;\; j=1, 2, \cdots, m.
	\end{align*}
	Since $ x\in [0,1] $, a simple computation shows that
	\begin{align*}
		x(1+x)^2A^2(\gamma)&\leq 4\left(\frac{3}{8}\right)^2\\ x(1+x)^2(1-x^2)^2A^4(\gamma)\vspace{0.2in}&\leq 4\left(\frac{3}{8}\right)^4\\&\vdots\\ x(1+x)^2(1-x^2)^{2m-2}A^{2m}(\gamma)&\leq 4\left(\frac{3}{8}\right)^{2m}.
	\end{align*}
It is easy to see that
\begin{align*}
J^{\prime}(x)&=\frac{2}{(1+x)^2}\bigg(1-2c_1x(1+x)^2A^2(\gamma)-6c_2 x(1+x)^2(1-x^2)^2A^4(\gamma)\\&\quad\quad-\cdots-2(2m-1)c_mx(1+x)^2(1-x^2)^{2m-2}A^{2m}(\gamma)\bigg)\\&\geq  \frac{2}{(1+x)^2}\left(1-\left(8c_1\left(\frac{3}{8}\right)^2+24c_2\left(\frac{3}{8}\right)^4+\cdots+8(2m-1)c_m\left(\frac{3}{8}\right)^{2m}\right)\right)\\&\geq 0,\;\;\;\;\; \text{if}\;\; 8c_1\left(\frac{3}{8}\right)^2+24c_2\left(\frac{3}{8}\right)^4+\cdots+8(2m-1)c_m\left(\frac{3}{8}\right)^{2m}\leq 1.
\end{align*}
Therefore, $ J(x) $ is an increasing function in $ [0,1] $ for 
\begin{equation*}
	8c_1\left(\frac{3}{8}\right)^2+24c_2\left(\frac{3}{8}\right)^4+\cdots+8(2m-1)c_m\left(\frac{3}{8}\right)^{2m}\leq 1.
\end{equation*}
	Hence, $ J(x)\leq 0 $ for all $ x\in [0,1] $ and $ \gamma\in [0,1) $. This completes the proof.
\end{pf}

\begin{pf}[\bf Proof of Theorem \ref{th-3.3}]
	From \eqref{e-2.3}, we have
	\begin{equation}\label{e-4.12a}
		|a_n|\leq\lambda(1-|a_0|^2)\;\;\mbox{for}\;\; n\geq 1.
	\end{equation}
A simple calculation using \eqref{e-4.12a} shows that
\begin{align}\label{e-3.4a}	\frac{S_r}{\pi}&=\displaystyle\frac{1}{\pi}\iint\limits_{|z|<r}|f^{\prime}(z)|^2dxdy\leq\sum_{n=1}^{\infty}n|a_n|^2r^{2n}\\&\nonumber\leq \lambda^2(1-|a_0|^2)^2\sum_{n=1}^{\infty}nr^{2n}\\&\nonumber=\lambda^2(1-|a_0|^2)\left(\frac{r}{1-r^2}\right)^2.
\end{align}
\noindent For $ r\leq 1/(1+2\lambda) $, it is easy to see that
\begin{equation}\label{e-4.12b}
	\frac{r}{1-r^2}\leq\frac{1+2\lambda}{4(1+\lambda)\lambda}.
\end{equation}
In view of \eqref{e-4.12a}, \eqref{e-3.4a} and \eqref{e-4.12b}, we obtain
\begin{align*} 
	\sum_{n=0}^{\infty}|a_n|r^n+P\left(\frac{S_r}{\pi}\right)&= |a_0|+\sum_{n=1}^{\infty}|a_n|r^n+\sum_{j=1}^{m}\left(\frac{1+\lambda}{1+2\lambda}\right)^{2j}\left(\frac{S_r}{\pi}\right)^j\\&=1-\frac{(1-|a_0|^2)}{16^m}\left(\frac{16^m}{1+|a_0|}-\frac{16^m}{2}-\sum_{j=1}^{m}16^{m-j}(1-|a_0|^2)^{2j-1}\right)\\& :\;=1-\frac{(1-|a_0|^2)}{16^m} J_1(|a_0|),
\end{align*}
where 
\begin{equation*}
	J_1(x)=\frac{16^m}{1+x}-\frac{16^m}{2}-\sum_{j=1}^{m}16^{m-j}(1-x^2)^{2j-1}.
\end{equation*}
In particular, we have
\begin{equation*}
	J_1(0)=\frac{16^m}{2}-\frac{16^m-1}{15}=\frac{13}{30}16^m-1>0\;\;\mbox{and}\;\; J_1(1)=0.
\end{equation*}
A simple computation shows that 
\begin{align*}
	J_1^{\prime}(x)&=-\frac{16^m}{(1+x)^2}+2x\sum_{j=1}^{m}16^{m-j}(1-x^2)^{2j-2}\\&\leq-\frac{16^m}{4}+2\frac{16^m-1}{16-1}\\&=-\left(\frac{7}{60}16^m+\frac{2}{15}\right)<0.
\end{align*}
Therefore, $ J_1(x) $ is a decreasing function in $ [0,1] $ and hence $ J_1(x)\geq J_1(1)=0 $ for $ x\in[0,1]. $ Thus, we have 
\begin{equation*}
	\sum_{n=0}^{\infty}|a_n|r^n+P\left(\frac{S_r}{\pi}\right)\leq 1\;\; \mbox{for}\;\; r\leq\frac{1}{1+2\lambda}.
\end{equation*}
It is easy to see that the equality
\begin{equation*}
	\sum_{n=0}^{\infty}|a_n|\left(\frac{1}{1+2\lambda}\right)^n+P\left(\frac{S_{\frac{1}{1+2\lambda}}}{\pi}\right)=1
\end{equation*}
holds for the function $ f $ of the form $ f(z)=\sum_{n=0}^{\infty} a_nz^n$
if, and only if, $ f\equiv c $ with $ |c|=1. $ This completes the proof.
\end{pf}
\begin{pf}[\bf Proof of Theorem \ref{th-3.4}]
	From \eqref{e-2.3}, we know that
	\begin{equation*}
		|a_n|\leq\lambda(1-|a_0|^2)\;\;\mbox{for}\;\; n\geq 1.
	\end{equation*}
	Then for $ r\leq 1/(1+2\lambda) $, a simple computation shows that 
	\begin{align}\label{e-3.6a}
		|a_0|+\sum_{n=1}^{\infty}\bigg(|a_n|+\beta|a_n|^2\bigg)r^n & \leq |a_0|+\bigg(\lambda(1-|a_0|^2)+\beta\lambda^2(1-|a_0|^2)^2\bigg)\frac{r}{1-r}\\&\nonumber\leq |a_0|+\frac{1-|a_0|^2}{2}\left(1+\lambda\beta(1-|a_0|^2)\right)\\&\nonumber: = 1-\frac{1-|a_0|^2}{2}F_{1}(x),
	\end{align}
where $ F_{1}(x) $ is defined by
\begin{equation*}
	F_{1}(x)=\frac{2}{1+x}-1-\lambda\beta(1-x^2). 
\end{equation*}
Since $ 0\leq\beta\leq 1/\lambda $, it follows that $ F_{1}(0)=1-\lambda\beta\geq 0 $ and $ F_{1}(1)=0 $. 
Our aim is to show that $ F_{1}(x)\geq 0 $ for $ x\in [0,1] $. It is easy to see that
\begin{equation*}
	F^{\prime}_{1}(x)=-\frac{2}{(1+x)^2}+2\beta\lambda x\leq -2+2\beta\lambda\leq 0\;\; \mbox{for}\;\; 0\leq \beta\leq\frac{1}{\lambda}.
\end{equation*}
Therefore, $ F_{1}(x) $ is a decreasing function on $ [0,1] $ and hence we have $ F_{1}(x)\geq F_{1}(1)=0 $ for $ x\in [0,1] $.
In view of \eqref{e-3.6a}, we obtain
\begin{equation*}
	|a_0|+\sum_{n=1}^{\infty}\bigg(|a_n|+\beta|a_n|^2\bigg)r^n\leq 1\;\;\mbox{for}\;\; r\leq\frac{1}{1+2\lambda}.
\end{equation*}
It is easy to see that the equality \begin{equation*}
	|a_0|+\sum_{n=1}^{\infty}\bigg(|a_n|+\beta|a_n|^2\bigg)\left(\frac{1}{1+2\lambda}\right)^n=1
\end{equation*}
holds for the function $ f $ of the form $ f(z)=\sum_{n=0}^{\infty} a_nz^n$
if, and only if, $ f\equiv c $ with $ |c|=1. $ This completes the proof.
\end{pf}

\begin{pf}[\bf Proof of Theorem \ref{th-3.5}]
		From \eqref{e-2.3}, we have 
		\begin{equation*}
		|a_n|\leq\lambda(1-|a_0|^2)\;\;\mbox{for}\;\; n\geq 1.
	\end{equation*}
Then for $ r\leq 1/(1+2\lambda) $, a computation shows that
\begin{align*} &
	B_1(r):=\sum_{n=0}^{\infty}|a_n|r^n+\left(\frac{1+\lambda}{2\lambda(1+|a_0|)}+\frac{2}{3}\frac{(1+\lambda) r}{1-r}\right)\sum_{n=1}^{\infty}|a_n|^2r^{2n}\\&\leq |a_0|+\lambda(1-|a_0|^2)\frac{r}{1-r}+\left(\frac{1+\lambda}{2\lambda(1+|a_0|)}+\frac{2}{3}\frac{(1+\lambda) r}{1-r}\right)\lambda^2(1-|a_0|^2)^2\frac{r^2}{1-r^2}\\&\leq |a_0|+\frac{1-|a_0|^2}{2}+\left(\frac{1+\lambda}{2\lambda(1+|a_0|)}+\frac{2}{3}\frac{(1+\lambda) r}{1-r}\right)\frac{\lambda}{4(1+\lambda)}(1-|a_0|)^2(1+|a_0|)^2\\&=1-\frac{1-|a_0|^2}{2}\left(1-\frac{1+|a_0|}{2}-\left(\frac{1}{4(1+|a_0|)}+\frac{1}{6}\right)(1-|a_0|)^2(1+|a_0|)\right)\\& =1-\frac{1-|a_0|^2}{2} F_2(|a_0|),
\end{align*}
where $ F_2(x) $ is defined by
\begin{equation*}
	F_2(x)=\frac{2}{1+x}-1-\frac{1}{2}\left(\frac{1}{2(1+x)}+\frac{1}{3}\right)(1-x)^2(1+x).
\end{equation*}\begin{figure}[!htb]
\begin{center}
	\includegraphics[width=0.50\linewidth]{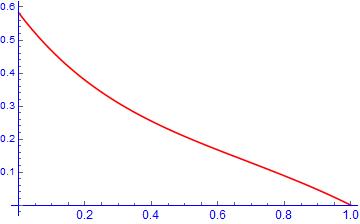}
\end{center}
\caption{Graph of $ F_2(x) $ in $ [0,1] $.}
\end{figure}
\begin{figure}[!htb]
	\begin{center}
		\includegraphics[width=0.50\linewidth]{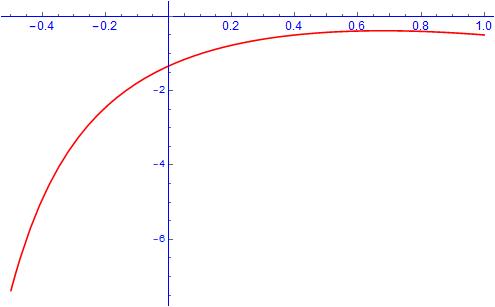}
	\end{center}
	\caption{The graph of the function $ F_2^{\prime}(x) $ in $ [0,1] $.}
\end{figure}
In order to prove that $ B_1(r)\leq 1 $, it is enough to show that $ F_2(x)\geq 0 $ for $ x\in [0,1] $. Clearly, $ F_2(1)=0 $ and
\begin{align*}
	 F_2(0)=1-\frac{1}{2}\left(\frac{1}{2(1+x)}+\frac{1}{3}\right)\geq\frac{7}{12}>0. 
\end{align*}
A simple calculation shows that for $ x\in[0,1] $,
\begin{equation*}
	F_2^{\prime}(x)=-\frac{2}{(1+x)^2}+\frac{(1-x)^2}{4(1+x)}-\frac{1}{2}\left(\frac{1}{2(1+x)}+\frac{1}{3}\right)(3x^2-2x-1)<0.
\end{equation*}
Therefore, the function $ F_2(x) $ is a decreasing function and hence we have $ F_2(x)\geq F_2(1)=0 $. As $  F_2(x)\geq 0 $ for $ x\in [0,1] $ and we obtain 
\begin{equation*}
	\sum_{n=0}^{\infty}|a_n|r^n+\left(\frac{1+\lambda}{2\lambda(1+|a_0|)}+\frac{2}{3}\frac{(1+\lambda) r}{1-r}\right)\sum_{n=1}^{\infty}|a_n|^2r^{2n}\leq 1\;\; \mbox{for}\;\; r\leq\frac{1}{1+2\lambda}.
\end{equation*}
Now the equality 
\begin{equation*}
	\sum_{n=0}^{\infty}|a_n|\left(\frac{1}{1+2\lambda}\right)^n+\left(\frac{1+\lambda}{2\lambda(1+|a_0|)}+\frac{2}{3}\frac{(1+\lambda)}{2\lambda}\right)\sum_{n=1}^{\infty}|a_n|^2\left(\frac{1}{1+2\lambda}\right)^{2n}=1
\end{equation*}
holds for a function $ f(z)=\sum_{n=0}^{\infty}a_nz^n $  in $ \mathcal{B}(\Omega) $ if, and only if, $ f\equiv c $ with $ |c|=1 $.
This completes the proof.
\end{pf}

\begin{pf}[\bf Proof of Theorem \ref{th-3.6}]
	Let $f \in \mathcal{B}(\Omega_{\gamma})$ and $ g(z)=f((z-\gamma)/(1-\gamma)) $. Then it is easy to see that $ g\in\mathcal{B}(\mathbb{D})$ and 
	\begin{align*}
		g(z)=\sum_{n=0}^{\infty}\frac{a_n}{(1-\gamma)^n}(z-\gamma)^n.
	\end{align*}
	Using Lemma \ref{lem-3.6}, we obtain 
	\begin{align*}
		\sum_{n=0}^{\infty}\frac{|a_n|}{(1-\gamma)^n}\rho^n+Q\left(\frac{S_{\rho}}{\pi}\right)\leq 1\;\; \text{for}\;\; \rho\leq\frac{1-\gamma^2}{3+\gamma}
	\end{align*} which is equivalent to 
	\begin{align}\label{e-4.15}
		\sum_{n=0}^{\infty}{|a_n|}\left(\frac{\rho}{(1-\gamma)}\right)^n+Q\left(\frac{S_{\rho}}{\pi}\right)\leq 1\;\; \text{for}\;\; \rho\leq\frac{1-\gamma^2}{3+\gamma}.
	\end{align}
	Set $ \rho=r(1-\gamma) $, then in view of \eqref{e-4.15}, we obtain 
	\begin{align}\label{e-3.8e}
		\sum_{n=0}^{\infty}{|a_n|}r^n+Q\left(\frac{S_{r(1-\gamma)}}{\pi}\right)\leq 1\;\; \text{for}\;\; r\leq\frac{1+\gamma}{3+\gamma}.
	\end{align} 
	To show the sharpness of the result, we consider the following function
	\begin{align*}
		f_a(z)=\frac{a-\gamma-(1-\gamma)z}{1-a\gamma-a(1-\gamma)}\;\;\text{for}\;\; z\in \Omega_{\gamma}\;\; \text{and}\;\; a\in (0,1).
	\end{align*}
	Define $ \phi_1 : \mathbb{D}\rightarrow\mathbb{D} $ by $ \phi_1(z)=(a-z)/(1-az) $ and $ \phi_2(z) :\Omega_{\gamma}\rightarrow\mathbb{D} $ by $ \phi_2(z)=(1-\gamma)z+\gamma $. 
	Then, the function $ f_a=\phi_1\circ \phi_2 $ maps $ \Omega_{\gamma} $ univalently onto $ \mathbb{D}. $ A simple computation shows that
	\begin{equation*}
		f_a(z)=\frac{a-\gamma-(1-\gamma)z}{1-a\gamma-a(1-\gamma)}=C_0-\sum_{n=1}^{\infty}C_nz^n\;\; \mbox{for}\;\; z\in\mathbb{D},
	\end{equation*} where $ a\in (0,1) $ and 
	\begin{equation*}
		C_0=\frac{a-\gamma}{1-a\gamma}\;\; \text{and}\;\; C_n=\frac{1-a^2}{a(1-a\gamma)}\left(\frac{a(1-\gamma)}{1-a\gamma}\right)^n.
	\end{equation*}
	A simple computation using \eqref{e-3.8e} shows that
	\begin{align*} &
		\sum_{n=0}^{\infty}{|a_n|}r^n+Q\left(\frac{S_{r(1-\gamma)}}{\pi}\right)\\&=\frac{a-\gamma}{1-a\gamma}+\left(\frac{1-a^2}{1-a\gamma}\right)\frac{(1-\gamma)r}{1-a\gamma-ar(1-\gamma)}+\frac{c_1r^2(1-a^2)^2(1-\gamma)^4}{((1-a\gamma)^2-a^2r^2(1-\gamma)^4)^2}\\[2mm]& \quad\quad+\frac{c_2r^4(1-a^2)^4(1-\gamma)^8}{((1-a\gamma)^2-a^2r^2(1-\gamma)^4)^4}+\cdots+\frac{c_mr^{2m}(1-a^2)^{2m}(1-\gamma)^{4m}}{((1-a\gamma)^2-a^2r^{2}(1-\gamma)^4)^{2m}} \\& :\;=1-(1-a)\Phi^{\gamma}_1(r),
	\end{align*}
	where 
	\begin{align*} 
		\Phi^{\gamma}_1(r)&=-\frac{(1+a)(1-\gamma)r}{(1-a\gamma-ar(1-\gamma))(1-a\gamma)}-\frac{c_1r^2(1-a)(1+a)^2(1-\gamma)^4}{((1-a\gamma)^2-a^2r^2(1-\gamma)^4)^2} \\[2mm]
		&\quad\quad -\frac{c_2r^4(1-a)^3(1+a)^4(1-\gamma)^8}{((1-a\gamma)^2-a^2r^2(1-\gamma)^4)^4}-\frac{c_mr^{2m}(1-a)^{2m-1}(1+a)^{2m}(1-\gamma)^{4m}}{((1-a\gamma)^2-a^2r^{2}(1-\gamma)^4)^{2m}}\\[2mm]&\quad\quad-\frac{1}{1-a}\left(\frac{a-\gamma}{1+a\gamma}-1\right).
	\end{align*}
	It is not difficult to show that $ \Phi^{\gamma}_1(r) $ is strictly decreasing function of $ r$ in $(0,1) $. Therefore,  for $ r>r_0=(1+\gamma)/(3+\gamma) $, we have $ \Phi^{\gamma}_1(r)<\Phi^{\gamma}_1(r_0) $. An elementary calculation shows that 
	\begin{align*}
		\lim_{a\rightarrow 1}\Phi^{\gamma}_1(r_0)=-\frac{2r_0}{(1-\gamma)(1-r_0)}+\frac{1+\gamma}{1-\gamma}=0.
	\end{align*}
	Therefore  $ \Phi^{\gamma}_1(r)< 0 $ for $ r>r_0 $. Hence $ 1-(1-a)\Phi^{\gamma}_1(r)>1 $ for $ r>r_0, $ which shows that $ r_0 $ is the best possible. 
\end{pf}
\begin{proof}[\bf Proof of Theorem \ref{th-3.1}]
	Let $ f=h+\overline{g} $ given by \eqref{e-2.4a} be harmonic in $ \Omega_{\gamma} $ with $ |h(z)|+|g(z)|\leq 1 $ in $ \Omega_{\gamma} $. Let 
	\begin{equation}\label{e-4.1}
		\psi_1(z)=h(z)+\epsilon g(z),
	\end{equation}
	where $ \epsilon\; (|\epsilon|=1) $ be arbitrary. Note that $ \psi_1 $ is analytic in $ \Omega_{\gamma} $ and $ |\psi_1(z)|\leq 1 $ in $ \Omega_{\gamma} $. From \eqref{e-4.1}, we have
	\begin{equation*}
		\psi_1(z)=a_0+\sum_{n=1}^{\infty}(a_n+b_n)z^n\quad\mbox{for}\quad z\in\mathbb{D}.
	\end{equation*}
	In view of Lemma \ref{lem-2.1}, we obtain 
	\begin{equation}
		\label{e-4.2} |a_n+\epsilon b_n|\leq\frac{1-|a_0|^2}{1+\gamma}\quad\mbox{for}\quad n\geq 1.
	\end{equation}
	Since $ \epsilon\; (|\epsilon|=1) $ is arbitrary, we have
	\begin{equation}\label{e-4.3}
		|a_n+b_n|\leq\frac{1-|a_0|^2}{1+\gamma}\quad\mbox{for}\quad n\geq 1.
	\end{equation}
	A simple computation using \eqref{e-4.3} shows that 
	\begin{equation}\label{e-4.4}
		|a_0|+\sum_{n=1}^{\infty}\bigg(|a_n|+|b_n|\bigg)r^n\leq |a_0|+\frac{1-|a_0|^2}{1+\gamma}\sum_{n=1}^{\infty}r^n=|a_0|+\frac{1-|a_0|^2}{1+\gamma}\left(\frac{r}{1-r}\right).
	\end{equation}
	Therefore, the inequality \eqref{e-4.4} is less than or equal to $ 1 $, provided 
	\begin{equation*}
		|a_0|+\frac{1-|a_0|^2}{1+\gamma}\left(\frac{r}{1-r}\right)\leq 1\quad\mbox{for}\quad r\leq r_0:=\frac{1+\gamma}{3+\gamma}.
	\end{equation*}
	Hence, we have 
	\begin{equation*}
		|a_0|+\sum_{n=1}^{\infty}\bigg(|a_n|+|b_n|\bigg)r^n\leq 1\quad\mbox{for}\quad r\leq r_0:=\frac{1+\gamma}{3+\gamma}.
	\end{equation*}
	This completes the proof.
\end{proof}	
\begin{proof}[\bf Proof of Theorem \ref{th-3.2}]
	Let $ f=h+\overline{g} $ given by \eqref{e-2.4a} be harmonic in $ \Omega_{\gamma} $ with $ |h(z)|+|g(z)|\leq 1 $ in $ \Omega_{\gamma} $. Consider the function 
	\begin{equation}\label{e-4.5}
		\psi_2(z)=h(z)+\epsilon g(z),
	\end{equation}
	where $ \epsilon\; (|\epsilon|=1) $ is arbitrary. Note that $ \psi_2 $ is analytic in $ \Omega_{\gamma} $ and $ |\psi_2|\leq 1 $ in $ \Omega_{\gamma} $. From the equation \eqref{e-4.1}, we have
	\begin{equation*}
		\psi_2(z)=a_0+\sum_{n=1}^{\infty}(a_n+b_n)z^n\quad\mbox{for}\quad z\in\mathbb{D}.
	\end{equation*}
	In view of Lemma \ref{lem-2.1}, we obtain 
	\begin{equation}
		\label{e-4.6} |a_n+\epsilon b_n|\leq\frac{1-|a_0+\epsilon b_0|^2}{1+\gamma}\quad\mbox{for}\quad n\geq 1.
	\end{equation}
	Since $ \epsilon\; (|\epsilon|=1) $ is arbitrary, we obtain
	\begin{equation}\label{e-4.7}
		|a_n-\epsilon b_n|\leq\frac{1-|a_0-\epsilon b_0|^2}{1+\gamma}\quad\mbox{for}\quad n\geq 1.
	\end{equation}
	From \eqref{e-4.6} and \eqref{e-4.7}, we obtain
	\begin{equation}\label{e-4.8}
		|a_n+\epsilon b_n|+|a_n-\epsilon b_n|\leq\frac{2-|a_0+\epsilon b_0|^2-|a_0-\epsilon b_0|^2}{1+\gamma}\;\; \mbox{for}\;\; n\geq 1.
	\end{equation}
	It is easy to see that \eqref{e-4.8} is equivalent to
	\begin{equation}\label{e-4.9}
		|a_n+\epsilon b_n|+|a_n-\epsilon b_n|\leq\frac{2\left(1-|a_0|^2-|b_0|^2\right)}{1+\gamma}\;\; \mbox{for}\;\; n\geq 1.
	\end{equation}
	From \eqref{e-4.9}, we obtain
	\begin{equation*}
		\bigg(|a_n+\epsilon b_n|+|a_n-\epsilon b_n|\bigg)^2\leq\frac{4\bigg(1-|a_0|^2-|b_0|^2\bigg)^2}{(1+\gamma)^2}
	\end{equation*} 
	which is equivalent to
	\begin{equation}\label{e-3.18a}
		|a_n|^2+|b_n|^2+|a^2-e^{i2\theta}b_n^2|\leq \frac{2\bigg(1-|a_0|^2-|b_0|^2\bigg)^2}{(1+\gamma)^2}.
	\end{equation} 
	By choosing $ \theta=\pi/2+\arg(a_n)-\arg(b_n) $, it follows from \eqref{e-3.18a} that
	\begin{equation*}
		|a_n|^2+|b_n|^2\leq\frac{\bigg(1-|a_0|^2-|b_0|^2\bigg)^2}{\left(1+\gamma\right)^2}
	\end{equation*}
	and hence
	\begin{equation}\label{e-4.10}
		\sqrt{|a_n|^2+|b_n|^2}\leq\frac{1-|a_0|^2-|b_0|^2}{1+\gamma}\quad\mbox{for}\quad n\geq 1.
	\end{equation}
	A simple computation using \eqref{e-4.10} shows that
	\begin{align}
		\label{e-4.11} 
		\sqrt{|a_0|^2+|b_0|^2}&+\sum_{n=1}^{\infty}\bigg(\sqrt{|a_n|^2+|a_n|^2}\bigg)r^n\\&\leq |a_0|+\frac{1-|a_0|^2}{1+\gamma}\sum_{n=1}^{\infty}r^n\nonumber\\&\nonumber=|a_0|+\frac{1-|a_0|^2}{1+\gamma}\left(\frac{r}{1-r}\right).
	\end{align}
	Therefore, from \eqref{e-4.10} and \eqref{e-4.11}, we obtain
	\begin{align*}
		\sum_{n=0}^{\infty}(\sqrt{|a_n|^2+|b_n|^2})r^n&\leq\sqrt{|a_0|^2+|b_0|^2}+\frac{1-|a_0|^2-|b_0|^2}{1+\gamma}\left(\frac{r}{1-r}\right)\\&\leq  \sqrt{|a_0|^2+|b_0|^2}+\frac{1-\sqrt{|a_0|^2+|b_0|^2}}{1+\gamma}\left(\frac{2r}{1-r}\right)\\&\leq 1,
	\end{align*} 
	if $ 2r/(1+\gamma)(1-r)\leq 1 $, that is if $ r\leq r_0=(1+\gamma)/(3+\gamma). $ Therefore, we have
	\begin{equation*}
		|a_0|+\sum_{n=1}^{\infty}\bigg(|a_n|+|b_n|\bigg)r^n\leq 1\quad\mbox{for}\quad r\leq r_0=\frac{1+\gamma}{3+\gamma}.
	\end{equation*}
	Hence, we have
	\begin{equation*}
		\sqrt{|a_0|^2+|b_0|^2}+\sum_{n=1}^{\infty}\sqrt{|a_n|^2+|b_n|^2}r^n\leq 1 \quad\mbox{for}\quad |z|=r\leq r_0= \frac{1+\gamma}{3+\gamma}.
	\end{equation*}
	This completes the proof.
\end{proof}

\noindent\textbf{Acknowledgment:}  The first author is supported by the Institute Post Doctoral Fellowship of IIT Bhubaneswar, India, the second author is supported by SERB-MATRICS, and third author is supported by CSIR, India.

\end{document}